\documentclass[10pt]{amsart}
\usepackage{amssymb, amsmath}
\usepackage[mathscr]{eucal}
\usepackage{graphicx}

\begin{document}

\title[Darboux's rotation and translation fields]{New manifestations of the Darboux's rotation\\ and translation fields of a surface}
\author[Victor Alexandrov]{Victor Alexandrov}
\address{
Sobolev Institute of Mathematics, Koptyug ave., 4, Novosibirsk, 630090, Russia and
Department of Physics, Novosibirsk State University, Pirogov str., 2, Novosibirsk, 630090, Russia}
\email{alex@math.nsc.ru}
\thanks{The author was supported in part by the Russian State 
Program for Leading Scientific Schools, Grant~NSh--8526.2008.1
and Federal Program ``Research and educational 
resouráes of innovative Russia in 2009--2013,'' 
contract No.~02.740.11.0457.}
\date{October 26, 2009}
\begin{abstract}
We show how the rotation and translation fields of a surface, introduced by
G.~ Darboux, may be used to obtain short proofs of a well-known theorem
(that reads that the total mean curvature of a surface is stationary under an infinitesimal bending)
and a new theorem (that reads that every infinitesimal flex of any simply connected closed
surface is orthogonal to the surface at least at two points).
\end{abstract}
\maketitle

A vector field
$\boldsymbol{z}=\boldsymbol{z}(u,v)$  
defined on a smooth surface 
$S\subset \Bbb R^3$
with a position vector 
$\boldsymbol{r}=\boldsymbol{r}(u,v)$ 
is called an \textit{infinitesimal bending} of~$S$ provided that
$$\boldsymbol{r}_u\cdot\boldsymbol{z}_u=\bold 0,\quad
\boldsymbol{r}_u\cdot\boldsymbol{z}_v+\boldsymbol{r}_v\cdot\boldsymbol{z}_u
=\bold0,\quad
\boldsymbol{r}_v\cdot\boldsymbol{z}_v=\bold 0,\eqno(1)
$$
where 
$\cdot$ 
stands for the inner product of 
$\Bbb R^3$ 
and
$\boldsymbol{r}_u, \boldsymbol{r}_v, \boldsymbol{z}_u, \boldsymbol{z}_v$ 
stand for the partial derivatives of the fields 
$\boldsymbol{r}$ and $\boldsymbol{z}$
with respect to the local parameters $u$ and $v$ on~$S$.

A vector field $\boldsymbol{z}$
is called a {\it nontrivial} infinitesimal bending of~$S$
provided that $\boldsymbol{z}$ is not a field of velocity vectors corresponding to a rigid motin of~$S$.

A surface is called \textit{nonrigid} if it admits a nontrivial
infinitesimal bending.

Note that nonrigid compact  boundary-free surfaces in~
$\Bbb R^3$ 
do exist and were studied by many authors
(see,  e.\,g.,  \cite{Re62, Sa92, Tr80} and references given there).

Any vector field $\boldsymbol{z}$ on a smooth surface~$S$
with the position vector $\boldsymbol{r}$ generates a mapping
$\boldsymbol{\psi}: S\times (-1,1)\to\Bbb R^3$
defined by the formula 
$\boldsymbol{\psi}(\boldsymbol{r},t) = \boldsymbol{r}+t\boldsymbol{z}(\boldsymbol{r})$.
It is an easy exercise that the following two statements are equivalent to each other:

(i) $\boldsymbol{z}$ is an infinitesimal bending of ~$S$;

(ii) for every smooth curve
$\gamma\subset S$,
the variation of the length of~$\gamma$ (i.\,e., the derivative with respect to~$t$ 
of the length of the curve
$\{\boldsymbol{\psi}(\boldsymbol{r},t)\mid \boldsymbol{r}\in\gamma\}$) 
vanishes at
$t=0$.

Our study is based on the notions of the  rotation and 
translation fields of a surface introduced by
G.~ Darboux \cite{Da96}, see also \cite{Bl21, Ef48, Sa65}.
The construction of the rotation field is based on the following theorem.

\textbf{Theorem 1.}
\textit{Let $\boldsymbol{z}=\boldsymbol{z}(u,v)=\bigl(z_1(u,v),z_2(u,v),z_3(u,v)\bigr)$  
be an infinitesimal bending of a smooth surface $S\subset\Bbb R^3$ with a position vector
$\boldsymbol{r}=\boldsymbol{r}(u,v)=\bigl(r_1(u,v),r_2(u,v),r_3(u,v)\bigr)$.
Then there exists a uniquely determined vector field 
$\boldsymbol{y}=\boldsymbol{y}(u,v)=\bigl(y_1(u,v),y_2(u,v),y_3(u,v)\bigr)$
on~$S$ such that $d\boldsymbol{z}=\boldsymbol{y}\times d\boldsymbol{r}$,
where $\times$ stands for the cross product on $\Bbb R^3$ or, equivalently,
such that $\boldsymbol{z}_u=\boldsymbol{y}\times\boldsymbol{r}_u$ and
$\boldsymbol{z}_v=\boldsymbol{y}\times \boldsymbol{r}_v$.}

\textit{Proof.} 
Rewriting the equation $\boldsymbol{z}_u=\boldsymbol{y}\times\boldsymbol{r}_u$
in the coordinate form
$$
\begin{array}{ll}
\boldsymbol{y}\times\boldsymbol{r}_u &=
\begin{pmatrix}
\boldsymbol{i} & \boldsymbol{j} & \boldsymbol{k} \\
y_1            & y_2            & y_3            \\
r_{1u}         & r_{2u}         & r_{3u}
\end{pmatrix}\\
 &=
\bigl(
y_2r_{3u}-y_3r_{2u}, y_3r_{1u}-y_1r_{3u}, y_1r_{2u}-y_2r_{1u}
\bigr)=
(z_{1u},z_{2u},z_{3u})=\boldsymbol{z}_u,
\end{array}
$$
we may treat it as the following system of linear algebraic equations
$$
\begin{pmatrix}
0 & r_{3u} & -r_{2u}\\
-r_{3u} & 0 & r_{1u}\\
r_{2u} & -r_{1u} & 0
\end{pmatrix}
\begin{pmatrix}
y_1 \\
y_2 \\
y_3
\end{pmatrix}
=
\begin{pmatrix}
z_{1u} \\
z_{2u} \\
z_{3u}
\end{pmatrix}.\eqno(2)
$$ 

The determinant of the $3\times 3$ matrix in (2) is equal to zero.
Hence, for some right-hand side vectors, (2) has no solutions while
for the others it has more than one solution.
We find conditions for the solvability of (2) using the Fredholm alternative.
In fact, it is easy to check that the corresponding homogeneous adjacent system 
$$
\begin{pmatrix}
0 & -r_{3u} & r_{2u}\\
r_{3u} & 0 & -r_{1u}\\
-r_{2u} & r_{1u} & 0
\end{pmatrix}
\begin{pmatrix}
p_1 \\
p_2 \\
p_3
\end{pmatrix}
=
\begin{pmatrix}
0 \\
0 \\
0
\end{pmatrix}
$$ 
has only one linearly indpendent solution, e.\,g.,
$(p_1,p_2,p_3)=(r_{1u}, r_{2u}, r_{3u})=\boldsymbol{r}_u$.
Consequently, (2) has a solution if and only if $\boldsymbol{r}_u\cdot\boldsymbol{z}_u=\boldsymbol{0}$.
This condition is fullfilled because of (1).
Moreover, any solution $\boldsymbol{y}$ to (2) is given by the formula
$\boldsymbol{y}=\boldsymbol{y}_1 +C_1\boldsymbol{r}_u$,
where $\boldsymbol{y}_1$ is any particular (or ``fixed'') solution to (2) 
and $C_1$ is an arbitrary constant.

Similarly, the equation $\boldsymbol{z}_v=\boldsymbol{y}\times \boldsymbol{r}_v$
is solvable because $\boldsymbol{r}_v\cdot\boldsymbol{z}_v=\boldsymbol{0}$
and its every solution $\boldsymbol{y}$ is given by the formula
$\boldsymbol{y}=\boldsymbol{y}_2 +C_2\boldsymbol{r}_v$,
where $\boldsymbol{y}_2$ is its any particular solution 
and $C_2$ is an arbitrary constant.

Since the vectors $\boldsymbol{r}_u$ and $\boldsymbol{r}_v$ are not collinear,
it follows that the system of equations 
$\boldsymbol{z}_u=\boldsymbol{y}\times\boldsymbol{r}_u$ and
$\boldsymbol{z}_v=\boldsymbol{y}\times \boldsymbol{r}_v$ has a unique solution $\boldsymbol{y}$.
\hfill$\square$

\textit{Remark.} From the analytical point of view the vector field $\boldsymbol{y}$ is more 
convenient than $\boldsymbol{z}$; this follows from the fact that the following statements are
equivalent to each other \cite{Co36}:

(a) an infinitesimal bending $\boldsymbol{z}$ is trivial;

(b) there are constant vector fields $\boldsymbol{a}$ and $\boldsymbol{b}$ such that
$\boldsymbol{z}=\boldsymbol{a}+\boldsymbol{b}\times\boldsymbol{r}$;

(c) $\boldsymbol{y}$ is a constant vector field.

\textbf{Definition.}
The vector field  $\boldsymbol{y}$, whose existence is istablished in Theorem 1,
is called the \textit{rotation field} of the surface~$S$ under the infinitesimal bending~ $\boldsymbol{z}$.
The vector field  $\boldsymbol{s}$, defined by the formula 
$\boldsymbol{s}=\boldsymbol{z}-\boldsymbol{y}\times\boldsymbol{r}$,
is called the \textit{translation field} of the surface~$S$ under the infinitesimal bending~ $\boldsymbol{z}$.

The rotation and translation fields were invented by G.~Darboux \cite{Da96}, who discovered 
a beautiful algebraic construction called ``the Darboux crown'', later studied by many authors,
see, e.\,g., \cite{Sa65, Sa92} and the references cited therein.
A briliant application of the rotation field were found by W.~Blaschke \cite{Bl21},
who proposed the simplest known proof of the rigidity of smooth ovaloids. Later Blaschke's proof
was popularised by many authors, see, e.\,g., \cite{Co36, Ef48}. 
The translation field played important role in the study by E. Rembs \cite{Re30} and 
R. Sauer \cite{Sa34}, who, among other things, have proved that a projective 
image of a nonrigid surface is a nonrigid surface again (this property attracts 
attantion of modern geometers too, see, e.\,g., \cite {Iz09}).
Among contemporary authors who use the rotation and
translation fields we can mention Ph.G. Ciarlet and O. Iosifescu \cite{CI09}.

The aim of the present paper is to show that the rotation and translation fields of 
a surface may be used to obtain short proofs of a well-known Theorem 3
(that reads that the total mean curvature is stationary under an infinitesimal bending)
and a new Theorem 4 (that reads that every infinitesimal bending of any simply connected closed
surface is orthogonal to the surface at least at two points). 

We study behaviour of the total mean curvature first.
If~
$S$ 
is an oriented surface in~
$\Bbb R^3$
then its {\it total mean curvature} is given by the classical formula
$$
H(S)=\int\limits_{S}
\frac{1}{2}\bigl(\kappa_1(\boldsymbol{r})+\kappa_2(\boldsymbol{r})\bigr)\,dS,\eqno(3)
$$
where
$\kappa_1(\boldsymbol{r})$ 
and 
$\kappa_2(\boldsymbol{r})$
are the principal curvatures of
$S$
at the surface point $\boldsymbol{r}$.
Similarly to (ii),
the {\it variation} $H'(S)$ of the total mean curvature of~ 
$S$ under the infinitesimal bending $\boldsymbol{z}$ is given by the formula
$$
H'(S)=\dfrac{d}{dt}\biggr|_{t=0}H\bigl(\boldsymbol{\psi}(S,t)\bigr),
\eqno(4)
$$
where
$\boldsymbol{\psi} (S,t)=\bigl\{\boldsymbol{\psi}(\boldsymbol{r},t)
=\boldsymbol{r}+t\boldsymbol{z}(\boldsymbol{r})\mid 
\boldsymbol{r}\in S, \ 0\leqslant t\leqslant 1\bigr\}$.

\textbf{Theorem 2.}
\textit{For every compact oriented smooth surface~
$S$ 
in~  
$\Bbb R^3$
and any its infinitesimal bending~ 
$\boldsymbol{z}$, 
the variation of the total mean curvature of~
$S$
is equal to the negative one half of the line integral of the rotation field~
$\boldsymbol{y}$ 
over the boundary~
$\partial S$
of
$S$, 
i.\,e.,
$$
H'(S)=-\frac12\int\limits_{\partial S}\boldsymbol{y}\cdot d\boldsymbol{r}.
$$}

Of course, here
$\partial S$
is oriented to be compatible with the orientation of~
$S$.

\textit{Proof.} 
It suffice to prove Theorem 2 ``locally,'' i.\,e., for~ 
$S$ 
covered by a single chart.
In particular, we may assume that~ 
$S$ 
is parameterized by 
$\boldsymbol{r}=\boldsymbol{r}(u,v)=\bigl(u,v,f(u,v)\bigr)$,
$(u,v)\in D\subset\Bbb R^2$. 
In agreement with standard notation, we put
$\boldsymbol{z}=(\xi,\eta,\zeta)$.
Then~(1) take the form
$$
\xi_u=-f_u\zeta_u,\quad
\xi_v+\eta_u=-f_v\zeta_u-f_u\zeta_v,\quad
\eta_v=-f_v\zeta_v,
\eqno(5)
$$
and the equations 
$\boldsymbol{z}_u=\boldsymbol{y}\times\boldsymbol{r}_u$,
$\boldsymbol{z}_v=\boldsymbol{y}\times \boldsymbol{r}_v$, 
defining the rotation field 
$\boldsymbol{y}$,
take the form
$$
f_uy_2=\xi_u,\quad
-f_uy_1+y_3=\eta_u,\quad
-y_2=\zeta_u,\eqno(6)
$$
$$
f_vy_2-y_3=\xi_v,\quad
-f_vy_1=\eta_v,\quad
y_1=\zeta_v.\eqno(7)
$$
Taking into account~(5), we find the followig solution to~(6) and (7):
$$
y_1=\zeta_v,\quad
y_2=-\zeta_u,\quad
y_3=\eta_u+f_u\zeta_v=-\xi_v-f_u\zeta_u.\eqno(8)
$$
Now, direct calculations show that
$$
\int\limits_{\partial S}\boldsymbol{y}\cdot d\boldsymbol{r}
=\int\limits_{\partial D}
(y_1+f_uy_3)\,du +(y_2+f_vy_3)\,dv.\eqno(9)
$$
Applying the Green's theorem
$$
\int\limits_{\partial D}P\,du+Q\,dv=\iint\limits_D\biggl
(\frac{\partial Q}{\partial u}-
\frac{\partial P}{\partial v}\biggr)\,dudv
$$
to the right-hand side integral in~(9), we obtain
$$
\int\limits_{\partial S}\boldsymbol{y}\cdot d\boldsymbol{r}=
\iint\limits_D
\bigl[
(y_{2u}+f_{uv}y_3+f_vy_{3u})-
(y_{1v}+f_{uv}y_3+f_uy_{3v})
\bigr]\,dudv.
$$
Using (8) and the formulas $\eta_{uu}=-f_{uu}\zeta_v-f_v\zeta_{uv}$,
$\xi_{vv}=-f_{vv}\zeta_u-f_u\zeta_{vv}$ (that may be obtained from (5) by means
of differentiation with respect to $u$ and $v$), we get
$$
\int\limits_{\partial S}\boldsymbol{y}\cdot d\boldsymbol{r}=
-\iint\limits_D
\bigl[
(1+f_v^2)\zeta_{uu}-2f_uf_v\zeta_{uv}+(1+f_u^2)\zeta_{vv}
\bigr]\,dudv.
$$
To conclude the proof, it remains to note that the right-hand side of the last formula is equal to
$-2H'(S)$, as it follows from (4) by straightforward calculations.
\hfill$\square$

\textit{Remark.}
In \cite{Al09} the reader may find another representation of the variation of the total
mean curvature in terms of a line integral.

\textbf{Theorem 3.}
\textit{For every compact oriented boundary-free smooth surface~ $S$ in 
$\Bbb R^3$
and any its infinitesimal bending, the variation of the total mean curvature of~ $S$ is equal to zero.}

\textit{Proof.}
Immediately follows  from Theorem 2.
\hfill$\square$

\textit{Remark.}
In fact, Theorem 3 was proven by other authors in a much more general situation, namely, for piecewise smooth hypersurfaces in multidimensional Euclidean and Lobachevskij spaces, see \cite{AR98, SS03, So04}.
But their proofs are much more complicated.

\textbf{Theorem 4.}
\textit{Let $\boldsymbol{r}$ be the position vector of a simply connected smooth surface $S\subset\Bbb R^3$,
let $\boldsymbol{s}$ be the translation field of~ $S$ under an infinitesimal bending $\boldsymbol{z}$ of~$S$,
and let $\Delta\subset S$ be a domain with smooth boundary and compact closure.
Then
$$
\int\limits_{\partial\Delta}\boldsymbol{s}\cdot d\boldsymbol{r}=\boldsymbol{0}.
$$}

\textit{Proof.}
It suffice to consider the case when~
$S$ 
is parameterized by 
$\boldsymbol{r}=\boldsymbol{r}(u,v)=\bigl(u,v,f(u,v)\bigr)$,
$(u,v)\in D\subset\Bbb R^2$
and $\boldsymbol{r}$ maps $D$ onto $\Delta$.
As usual, we put
$\boldsymbol{z}=(\xi,\eta,\zeta)$.

Using the definition
$\boldsymbol{s}=\boldsymbol{z}-\boldsymbol{y}\times\boldsymbol{r}$
and taking into account (8),
we easily find the coordinates of~$\boldsymbol{s}=(s_1,s_2,s_3)$:
$$
\begin{array}{ll}
 & s_1=\xi+f\zeta_u+v\eta_u+vf_u\zeta_u,\\
 & s_2=\eta+f\zeta_v-u\eta_u-uf_u\zeta_v,\\
 & s_3=\zeta-u\zeta_u-v\zeta_v.
\end{array}\eqno(10)
$$

We have
$$
\int\limits_{\partial\Delta}\boldsymbol{s}\cdot d\boldsymbol{r}
=\int\limits_{\partial D}
(s_1+f_us_3)\,du +(s_2+f_vs_3)\,dv.\eqno(11)
$$
Applying the Green's theorem
to the right-hand side integral in~(11),
using~(10),
and taking into account the formulas
$\xi_v+\eta_u=-f_v\zeta_u-f_u\zeta_v$,
$\eta_{uu}=-f_{uu}\zeta_v-f_v\zeta_{uv}$, and
$\eta_{uv}=-f_{uv}\zeta_v-f_v\zeta_{uv}$,
we obtain after simplifications that the function under the 
sign of the double integral vanishes identically.
Hence, the left-hand side in (11) is equal to zero.
\hfill$\square$

\textbf{Definition.} 
Let $\boldsymbol{r}$ be the position vector of a point of 
a connected, simply connected surface~$S$ in $\Bbb R^3$.
Let $\boldsymbol{s}$ be the translation field of~$S$
under an infinitesimal bending $\boldsymbol{z}$.
By definition, put
$$
\varphi(\boldsymbol{r})=
\int\limits_{\gamma}\boldsymbol{s}\cdot d\boldsymbol{r},\eqno(12)
$$
where $\gamma\subset S$ is any smooth curve with the start point 
$\boldsymbol{r}_0$ and finish point $\boldsymbol{r}$
($\boldsymbol{r}_0$ is supposed to be an arbitrary ``fixed'' point of~$S$).

\textit{Remark.}
It follows from Theorem 4 that the above definition is consistent,
i.\,e., as soon as $\boldsymbol{r}_0$ is fixed, the value
$\varphi(\boldsymbol{r})$ is defined correctly, in particular,
it does not depend on the choise of~$\gamma$.
Of course, the line integral in~(12) can be treated as work of the vector field
$\boldsymbol{s}$ along the curve~$\gamma$.
Thus, Theorem 4 reads that $\boldsymbol{s}$ is a potential field
and $\varphi$ is its potential function.

A simple consiquence of the existence of the potential function  $\varphi: S\to\mathbb R$
is given by the following theorem.

\textbf{Theorem 5.}
\textit{For every simply connected compact boundary-free smooth surface $S\subset\Bbb R^3$
and every its infinitesimal bending~ $\boldsymbol{z}$, there are at least two distinct
points where~ $\boldsymbol{z}$ is orthogonal to~ $S$.}

\textit{Proof.}
We prove that, on every connected component of $S$, the two points under study
are the points where the potential fuction $\varphi$ attains its maximum or minimum.

Let $A\in S$ be a point where~$\varphi$ attains a local extremum.
Choose a special coordunate system around $A$ such that
$\boldsymbol{r}=\boldsymbol{r}(u,v)=\bigl(u,v,f(u,v)\bigr)$,
$(u,v)\in D\subset\Bbb R^2$, $\boldsymbol{r}(0,0)=A$, and
$f(0,0)=f_u(0,0)=f_v(0,0)=0.$
Since
$$
\varphi(\boldsymbol{r})=
\int\boldsymbol{s}\cdot d\boldsymbol{r}=
\int(s_1+f_us_3)\,du +(s_2+f_vs_3)\,dv,
$$
it follows that 
$\varphi_u=s_1+f_us_3$ and $\varphi_v=s_2+f_vs_3$.
Now if we recall~(10), we get
$\varphi_u(0,0)=s_1(0,0)=\xi(0,0)$ and
$\varphi_v(0,0)=s_2(0,0)=\eta(0,0)$.
On the other hand, 
$\varphi_u(0,0)=\varphi_v(0,0)=0$.
Thus 
$\xi(0,0)=\eta(0,0)=0$ 
and the vector field~ 
$\boldsymbol{z}$
is orthogonal to~ $S$ at~ $A$.
\hfill$\square$

\textit{Remark.}
One may be tempted to consider Theorem 5 as a special case of theorems about zeros of Killing vectors fields (i.\,e., infinitesimal isometries of Riemannian manifolds) proved by S.~Kobayashi \cite{Ko72}. However this is not correct, because, in general, a field of infinitesimal bending of a surface $S\subset\Bbb R^3$ neither is a tangential vector field on~$S$ nor generates a Killing field on~$S$ in any natural way.

\bibliographystyle{plain}
\bibliography{bib_alex}

\end{document}